\documentclass[final,epsfig,amsfont]{article}
\usepackage{latexsym}
\usepackage[dvips]{graphicx}
\newcommand{\eeq}{\end{equation}}
\newcommand{\beq}{\begin{equation}}
\newcommand{\nuq}[1]{\label{#1} \eeq}

\newtheorem{definition}{Definition}
\newtheorem{lemma}{Lemma}
\newtheorem{corollary}{Corollary}
\newtheorem{problem}{Problem}
\newtheorem{hypothesis}{Hypothesis}
\newcommand{\bseq}{$\{b_j\}_{j \in \mathbf{N}}$}
\newtheorem{theorem}{Theorem}
\newtheorem{example}{Example}
\begin{document}
\title{
Numerical computation of the isospectral torus of finite gap sets and of IFS Cantor sets}
\author{
Giorgio Mantica \\
Center for Non-linear and Complex Systems, \\ Department of Physics and Mathematics, \\Universit\`a dell' Insubria, 22100 Como, Italy \\ and CNISM unit\`a di Como,
and  I.N.F.N. sezione di Milano\thanks{Computations for this paper were performed at the INFN computer center at Pisa.}.}

\date{}
\maketitle
\begin{abstract}
We describe a numerical procedure to compute the so--called isospectral torus of
finite gap sets, that is, the set of Jacobi matrices whose essential spectrum is composed of finitely many intervals. We also study numerically the convergence of specific Jacobi matrices to their isospectral limit.  We then extend the analysis to the definition and computation of an ``isospectral torus'' for Cantor sets in the family of Iterated Function Systems. This analysis is developed with the ultimate goal of attacking numerically the conjecture that the Jacobi matrices of I.F.S. measures supported on Cantor sets are asymptotically almost-periodic.\\
{\em keywords: Potential Theory \and Orthogonal Polynomials \and Iterated Function Systems \and Almost periodic Jacobi matrices.}\\
{\em MATH Subj. Class. 42C05; 31A15; 47B36; 81Q10; 30C30}

\end{abstract}
\vspace{1.0cm}
{\em Dedicated to Ed Saff on his seventieth}

\section{Introduction}
\label{sec1}

Harmonic analysis have always been one of Ed's interests, and while working on problems of theoretical sophistication, he never disdained the applicative side of the theory. Years ago we had the occasion to discuss the usage of Szeg\"o polynomials in the determination of harmonic frequencies embedded in a real (or complex) signal \cite{edfilt,edfilt2,filt1,filt2}. Ed has been a pioneer of the mathematical analysis of this approach, that is now being applied to challenging problems like gravitational waves detection by a long-time friend of Ed's, Daniel Bessis \cite{waves}.

In this paper I want to discuss a theoretical problem in which harmonic analysis plays a fundamental role, and where numerical experiments are important, for lack of a definitive theory. About twenty years ago \cite{physd1} I conjectured that the entries of the Jacobi matrix associated with the devil's staircase measure are almost periodic (all terms to be defined precisely later on), and I provided the first numerical evidence to that aim. The conjecture was rather natural, since another example of singular continuous measure supported on a Cantor set, that is, Julia set measures, possess almost periodic Jacobi matrices, provably so at least in a certain range of parameters \cite{belli2,belli,baker,pierre,mijeff}, and the so-called {\em almost periodic flu} \cite{last,testard,russell,cycon,testard} had just raged. Moreover, in the inverse direction, almost periodic Jacobi matrices, like those constructed according to the Fibonacci rule have a spectral measure supported on Cantor sets with a clear hierarchical structure \cite{fibo1}.

Although this conjecture had been reported repeatedly afterwards \cite{dubna,mobius,physd2}, orthogonal polynomials on I.F.S. Cantor sets remained confined to the study of quantum evolution of systems with singular spectra, see {\em e.g.} \cite{etna} in  Ed's {\em 60th} anniversary volume and references therein.
This topic, and the conjecture, were recently brought to evidence by two papers that were motivated on the one hand \cite{stric0} by the search for approximate identities with Dirichlet kernels, but went much further, and on the other hand \cite{simonlast} by the infinite generalization of results for the so--called finite gap sets, that is, set of intervals, that opened new perspectives in this problem.
It is therefore timely to return to this problem. Although concerned with this dated question, the present paper is not simply review, but it introduces new ideas and techniques, as the occasion of Ed's seventieth deserves. As they say, happy continuation, Ed!

\section{Jacobi matrices on IFS sets: definitions and summary of the paper}
\label{sec2}
Orthogonal polynomials,  $\{p_j(\mu;s)\}_{j \in {\bf N}}$, of a positive Borel measure $\mu$ supported on a compact subset $E$ of the real axis are defined in a straightforward way
by the relation $\int p_j(\mu;s) p_m(\mu;s) d \mu(s) = \delta_{jm}$, where $\delta_{jm}$ is the Kronecker delta. The well known three-terms recurrence relation
\begin{equation}
\label{nor2}
   s p_j (\mu;s) = b_{j+1} p_{j+1}(\mu;s)
   + a_j p_j(\mu;s)  + b_{j} p_{j-1}(\mu;s),
\end{equation}
initialized by $b_0=0$ and
$p_{-1}(\mu;s) = 0$, $p_0(\mu;s) = 1$, can be formally encoded
in the Jacobi matrix $J(\mu)$:
\beq
   J(\mu) :=
           \left(   \begin{array}{ccccc}  a_0  & b_1 &      &      &       \cr
                               b_1  & a_1 & b_2  &      &       \cr
                      &   \ddots    &    \ddots      & \ddots &\cr
                               \end{array} \right) .
 \nuq{jame}
For compact support $E$ the moment problem is determined \cite{ach}, and the matrix $J(\mu)$
is in one--to--one relation with the measure $\mu$.

In this paper we will study the Jacobi matrix of a variety of measures supported on set of intervals, defined by Iterated Function Systems, reviewed
in Sect. \ref{sec-ifs}. This construction produces a sequence of finite unions of intervals, $\{E_n\}_{n \in {\bf N}}$, converging to a Cantor set when $n$ tends to infinity.
Each $E_n$ is a compact set, that carries a unique equilibrium measure, $\nu_n$.
A numerical procedure to compute this measure is described in Sect. \ref{sec-equi}, following \cite{dolo}. This procedure yields the set of {\em fundamental harmonic frequencies} that are the basis of further developments.
In sect. \ref{isotorfin} a novel numerical method for the computation of the isospectral torus of the sets $E^n$ is described. The infinite Jacobi matrices so computed are then analyzed in their harmonic components in sect. \ref{harman}, by using a refined technique outlined in Appendix. Validity of this technique is independently verified in the successive sect. \ref{sec-harm}, where we describe convergence to the isospectral torus of the Jacobi matrix elements of two significant I.F.S. examples. It might be that this is the first time that this convergence has been explicitly exhibited in a numerical example.

The final two sections touch upon two challenging problems. First, in sect. \ref{sec-isoinfty}, we propose a definition of the isospectral torus of infinite gap sets--Cantor sets generated by I.F.S. that naturally leads to numerical computations. Finally, in the closing sect. \ref{sec-if} we propose a framework to assess the validity of the conjecture mentioned in the introduction, that is, we discuss whether the Jacobi matrix of an I.F.S. balanced measure---like the devil's staircase measure---is asymptotically almost periodic. We study to what extent its matrix elements can be approximated by the trigonometric polynomials obtained by the harmonic analysis developed in the preceding sections.

\section{Attractors of Iterated Function Systems}
\label{sec-ifs}

Let us first construct finite gap sets, converging to a Cantor set. To do this, we use the simplest version of Iterated Function Systems (IFS) \cite{papmor,hut,dem,ba2}.
These are collections of contractive maps $\phi_i : {\bf R} \rightarrow
{\bf R}$, $i = 1, \ldots, M$, for which there exists a unique
set ${\mathcal A}$, called the {\em attractor} of the IFS, that solves
the equation
 \begin{equation}
 \label{attra}
    {\mathcal A}=\bigcup_{j=1,\ldots ,M}\;\phi_j({\mathcal A}) :=  \Phi ({\cal A}).
 \end{equation}
The attractor ${\mathcal A}$ can also be seen the limit (in the Hausdorff metric) of the sequence $E^n := \Phi^{n}(E^{0})$, where $E^{0}$ is any non-empty compact set:
$
{\cal A} = \lim_{n \to \infty} \Phi^{n}(E^{0}).
$
We now choose a finite number of affine maps of the form:
\begin{equation}
\label{mappi}
    \phi_j (s) = \delta_j (s - \gamma_j) + \gamma_j,  \;\;  j = 1, \ldots, M ,
\end{equation}
where $\delta_j$ are real numbers between zero and one, called {\em contraction ratios}, and $\gamma_{j}$ are real constants, that geometrically correspond to the fixed points of the maps. Furthermore, we choose the set $E^0$ as the interval $E^0=[\min \gamma_{j},\max \gamma_j]$, and we also impose that $\phi_i(E^0) \cap \phi_j(E^0) = \emptyset$, for $i \neq j$. This gives rise to what is called a {\em fully disconnected} IFS, whose attractor  ${\cal A}$ is a Cantor set. It follows from standard theory that $E^0$ is the convex hull of ${\cal A}$ and that
$E^n$ is the union of $M^n$ disjoint, closed intervals, that we will denote as $E^n_i = [\alpha_i,\beta_i]$:
\begin{equation}
\label{eq-int1}
    E^n = \Phi^n(E^0) = \bigcup_{i=1}^{M^n} [\alpha_i,\beta_i].
\end{equation}

In accordance with the spectral analysis of periodic solids, we will call the intervals $E^n_i$ at r.h.s. of eq. (\ref{eq-int1}) {\em bands} at generation (or level) $n$. It appears from eq. (\ref{eq-int1}) that all bands $E_i^n$ at level $n$ can be obtained by applying the transformations $\phi_j$, $j=1,\ldots,M$, to the bands at generation $n-1$. As a consequence, the length of these intervals, as well as their sum, tend to zero geometrically and the attractor is a Cantor set of null Lebesgue measure. Complementary to bands, in the convex hull of $E^n$, are
the so--called {\em gaps}, for which we will use the notation
\begin{equation}
\label{eq-gap1}
   G^n := \bigcup_i G^n_i = \bigcup_{i=1}^{M^n-1} (\beta_i,\alpha_{i+1}).
\end{equation}

\section{Equilibrium measure for a set of intervals}
\label{sec-equi}
Crucial in the following is the role of the equilibrium measure for $E^n$, for an electrostatic charge, with a logarithmic law of repulsion. We now recall briefly the main facts of this theory \cite{ran0,ed}. Let $\sigma$ any Borel probability measure, supported on $E^n$. The potential $V(\sigma;z)$, generated by $\sigma$ at the point $z$ in $\bf C$, is
 \begin{equation}
 \label{pote1}
    V(\sigma;z) := - \int_{E^n} \log |z-s| \; d \sigma(s).
 \end{equation}
The electrostatic energy ${\cal E}(\sigma)$ of the distribution $\sigma$ is the integral of $V(\sigma;z)$:
 \begin{equation}
 \label{pote2}
  {\cal E}(\sigma) := \int_{E^n}  V(\sigma;u) \; d\sigma(u)  = - \int_{E^n} \int_{E^n} \log |u-s| \; d\sigma(s)d\sigma(u).
 \end{equation}
The {\em equilibrium measure} $\nu_{E^n}$ associated with any compact domain $E^n$ is the unique measure that minimizes the energy ${\cal E}(\sigma)$, when this latter is not identically infinite. This being easily proven in the case of a finite union of intervals, $\mbox{Cap}(E^n) := e^{-{\cal E}(\nu_{E^n})}$  defines the capacity of the set $E^n$, and the Green's function can be written as
\begin{equation}
 \label{pote8}
  g(E^n;z) = - V(\nu_E^n;z) - \log (\mbox{Cap}(E^n)) = - V(\nu_{E^n};z) + {\cal E}(\nu_{E^n}).
 \end{equation}

From a theoretical point of view, the solution of the equilibrium problem for a finite union of $N$ intervals $[\alpha_i,\beta_i]$, $i = 1,\ldots,N$, is well known. Let's follow {\em e.g.} \cite{widom}. Define the polynomial $Y(z)$,
 \begin{equation}
 \label{meas2}
    Y(z) = \prod_{i=1}^{N} (z-\alpha_i)(z-\beta_i),
 \end{equation}
and its square root, $\sqrt{Y(z)}$, as the one which takes real values for $z$ real and large. Also, let the real number $\zeta_i$ belong to the open interval $(\beta_i,\alpha_{i+1})$, for $i=1,\ldots,N-1$ ({\em i.e.} to the gap $G^n_i$). Define $Z(z)$ as the monic polynomial of degree $N-1$ with roots at all $\zeta_i$'s:
 \begin{equation}
 \label{meas3}
    Z(z) =  \prod_{i=1}^{N-1} (z-\zeta_i).
 \end{equation}
\begin{theorem}
There exists a unique set of values $\{\zeta_i, i = 1,\ldots,N-1\}$ that solve the set of coupled, non--linear equations
 \begin{equation}
 \label{meas5}
    \int_{b_i}^{a_{i+1}} \frac{Z(s)}{\sqrt{|Y(s)|}} \; ds = 0, \;\; i = 1,\ldots,N-1.
 \end{equation}
In terms of these values, one can write the equilibrium measure of the set $E^n$: for simplicity of notation, we denote it by $\nu^n := \nu_{E^n}$:
 \begin{equation}
 \label{meas6z}
     d \nu^n (s) = \frac{1}{\pi} \sum_{i=1}^N \chi_{[\alpha_i,\beta_i]} (s)\;  \frac{|Z(s)|}{\sqrt{|Y(s)|}} \; ds.
 \end{equation}
\end{theorem}
It is apparent from the previous equation that the measure $\nu^n$ is absolutely continuous with respect to the Lebesgue measure on $E^n$.
Moreover, since the compact sets $E^n$ of positive capacity tend to $\cal A$ when $n$ grows, standard theory (see {\em e.g.} \cite{andriev}) assures that
\begin{theorem}
As $n$ tends to infinity,  $\nu_n$ tends weakly to $\nu_{\cal A}$, the equilibrium measure on the Cantor set $\cal A$.
\end{theorem}

In two previous works \cite{dolo,arx}, we have shown how to solve numerically the set of equations (\ref{meas5}) in a stable way, for sets of intervals of large cardinality.   In \cite{arx} we have also  quantified the rate of this convergence in terms of the orthogonal polynomials $p_j(\nu^n;z)$.

Of particular relevance are the integrals of $\nu^n$ over the intervals composing $E^n$, $\nu^n([\alpha_i,\beta_i])$.
When these integrals are rational numbers, Peherestorfer \cite{franz} has proven that there exists a strict-T polynomial on $E^n$ , that is, a polynomial with oscillation properties mimicking, and extending, those of the classical Chebyshev polynomials on a single interval. This is notably the case when the I.F.S. maps defining $E^n$ are the inverse branches of a polynomial dynamics with  a real Julia set (see \cite{belli,baker,pierre,mijeff,jeff2} and \cite{arx} for a numerical investigation). On the contrary, this is {\em not} the typical case of intervals generated by affine maps of the kind (\ref{mappi}). This fact has important consequences in the following. In fact, the harmonic analysis that we will develop features the {\em fundamental harmonic frequencies}:
\begin{equation}
 \label{meas8}
     \omega^n_{i} := \frac{1}{\pi} \int_{-\infty}^{\xi_i}  \frac{|Z(s)|}{\sqrt{|Y(s)|}} ds,
 \end{equation}
where $\xi_i$ is any point in the $i-th$ gap $[\alpha_i,\beta_{i+1}]$. Typically, for IFS affine maps, the values $\omega^n_i$, $i=1,\ldots,N-1$ are {\em not} rationally related.

\section{Set of Intervals -- The Isospectral torus}
\label{isotorfin}

The isospectral torus of Jacobi matrices associated with measures supported on a finite set of intervals arise naturally in the study of the asymptotics of orthogonal polynomials supported on these sets \cite{widom,nuttall,sasha1,franz,franzimrn,franz0,chris,chrisb,chris2,chris3}.  Therefore it has been extensively studied and many equivalent definitions shed light on its nature \cite{chris2b,simonlast}.
We want to adopt here a constructive definition \cite{alphonse},\cite{alphonse2} that reads as follows.
Let $Y(z)$ be defined as above, let $Z(\xi;z) =  \prod_{i=1}^{N-1} (z-\xi_i)$, where now $\xi_i$ are arbitrary values in the gaps and let $X$ be a polynomial constructed in such a way that the function
\begin{equation}
 \label{lupa1}
  M(z) = \frac{X(z)-\sqrt{Y(z)}}{Z(z)}
  \end{equation}
behaves as $1/z$ for large $z$ and $X(\xi_i) = \sigma_i \sqrt{Y(\xi_i)}$, with $|\sigma_i|=1$. Clearly, the $n-1$ numbers $\xi_i$ and the signs $\sigma_i$ label points in an $n-1$ dimensional torus. It can be shown that
\begin{theorem}
$M(z)$ is the Markov function (a.k.a. Stieltjes transform) of a positive measure $\theta$,
\begin{equation}
 \label{lupa2}
  M(z) =   \int \frac{ d \theta(s)}{z-s},
  \end{equation}
where $\theta$ is given by the sum of an absolutely continuous measure and of a linear combination of atomic measures supported at the zeros of the function $Z$:
\begin{equation}
 \label{meas6}
     d \theta (s) = \frac{A}{\pi} \sum_{i=1}^N \chi_{[\alpha_i,\beta_i]} (s)\;  \frac{\sqrt{|Y(s)|}}{|Z(\xi;s)|} \; ds + \sum_{i=1}^N B_i \delta_{\xi_i} .
 \end{equation}
\end{theorem}
Observe that the constants $A$ and $B_i$ depend on the parameters $\xi_i$ and $\sigma_i$ via eq. (\ref{lupa2}).
The absolutely continuous part contains a generalized Chebyshev measure of the second kind.
\begin{definition}
The isospectral torus of a finite set of intervals is the set of Jacobi matrices associated with the measures $\theta$ in eqs. (\ref{lupa1},\ref{lupa2},\ref{meas6}).
\label{defi-isofinite}
\end{definition}

This theorem prompts for a numerical computation of such Jacobi matrices, adapting the techniques of \cite{arx}. In essence, the theory run as follows:
\begin{itemize}
\item[$\bullet$] The absolutely continuous component in (\ref{meas6}) is written as the (weak) limit of a sequence of atomic measures (Lemma 2 in \cite{arx}), constructed using the properties of Chebyshev measures, of first kind for the equilibrium measure, of second kind for the present application.
    \item[$\bullet$] The (truncated) Jacobi matrix of a finite sum of atomic measures is computed adapting and enhancing existent procedures. These techniques allow for a uniform control of convergence to the desired matrix elements: we can estimate  $\max \{|a_j-\bar{a}_j|,|b_j-\bar{b}_j|, j=0,\ldots,J \}$, where the bar indicate the numerically computed matrix elements, and reduce this distance within acceptable limits, even for large values of $J$.
\end{itemize}
While our techniques are fully general, we need to choose a specific example for the numerical computations of this paper. At difference with Ref. \cite{arx}, where the ternary Cantor set has been chosen, we focus here on a different I.F.S., that does {\em not} enjoy of the useful symmetry properties of the former, and represents a more significant instance of the general case. This is the following
\begin{example} Let $M=2$, $\delta_1=0.34$, $\delta_2=0.52$, $\gamma_1=-1$, $\gamma_2=1$.
\label{examp-cantor}
\end{example}

\begin{figure}
\centerline{\includegraphics[width=.6\textwidth, angle = -90]{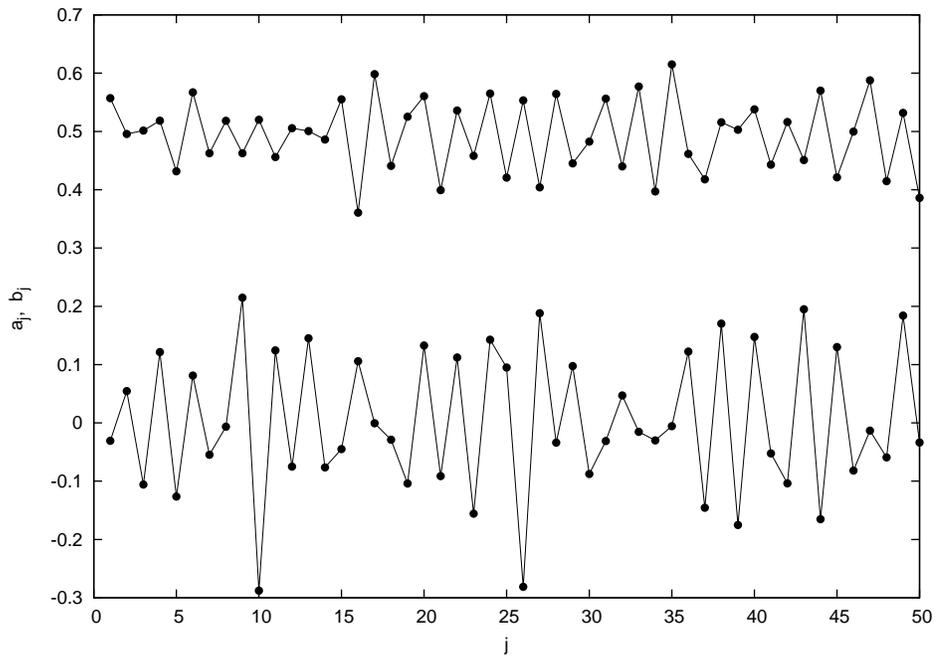}}
\caption{Diagonal (lower data) and outdiagonal (upper data) entries of a  Jacobi matrix in the isospectral torus of $E^4$ (IFS of Example \ref{examp-cantor}) .}
\label{fig-jaco1}
\end{figure}

\begin{figure}
\centerline{\includegraphics[width=.6\textwidth, angle = -90]{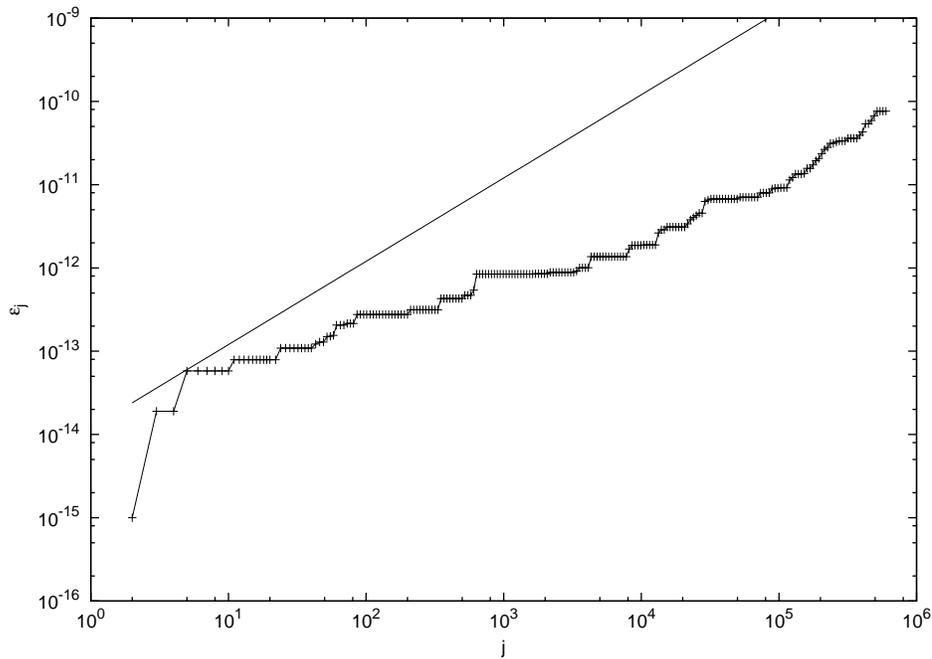}}
\caption{Error in the outdiagonal elements of the numerically computed Jacobi matrix of fig. \ref{fig-jaco1} versus $j$. The straight line has unit slope.}
\label{fig-jaco2}
\end{figure}

In figure \ref{fig-jaco1} we draw the initial parts of the sequences $\{a_j\}_{j=1}^\infty$, $\{b_j\}_{j=1}^\infty$ for a point in the isospectral torus of the set $E^4$ generated by the IFS of example \ref{examp-cantor}. Our technique allows for large--scale computations: in figure \ref{fig-jaco2} an estimate of the infinity distance between the numerically computed elements and the true ones, {\em i.e.} $\epsilon_j = \max \{|b_l-\bar{b}_l|, l \leq j\}$ is plotted versus $j$. We observe a less than linear growth of the error that permits to obtain extremely accurate values even for matrix sizes in the order of hundreds of thousands.

\section{Harmonic analysis of the isospectral torus of finite gap sets}
\label{harman}

In the previous section we have shown how to compute numerically the entries of Jacobi matrices in the isospectral torus of a finite number of intervals. It must be remarked that the same objects can also be expressed in terms of theta functions \cite{sasha1,sasha}, but this procedure seems to be numerically hard even for a small number of intervals. To the contrary, the approach adopted herein is capable of dealing satisfactorily with as many as $2^{18}$ intervals with moderate computer resources \cite{arx}. Before proceeding further with numerical techniques, we need an analytic result:

\begin{theorem}[\cite{franz0}]
There exist real analytic functions $\Phi$, $\Psi$ on the torus $[0,1]^{N-1}$ such that the recurrence coefficients $a_j, b_j$ of Jacobi matrices in the isospectral torus can be expressed as
\begin{equation}
\label{eq:phi1}
  a_j = \Phi (j \mathbf{\omega} + \mathbf{\phi}), \quad b_j = \Psi (j \mathbf{\omega} + \mathbf{\phi}),
\end{equation}
where $\mathbf{\omega}$ is the $n-1$ dimensional vector with components $\omega^n_i$ in eq. (\ref{meas8}), and $\mathbf{\phi}$ and $\mathbf{\psi}$ are shift vectors that depend on the particular matrix.
\end{theorem}

\begin{figure}
\centerline{\includegraphics[width=.6\textwidth, angle = -90]{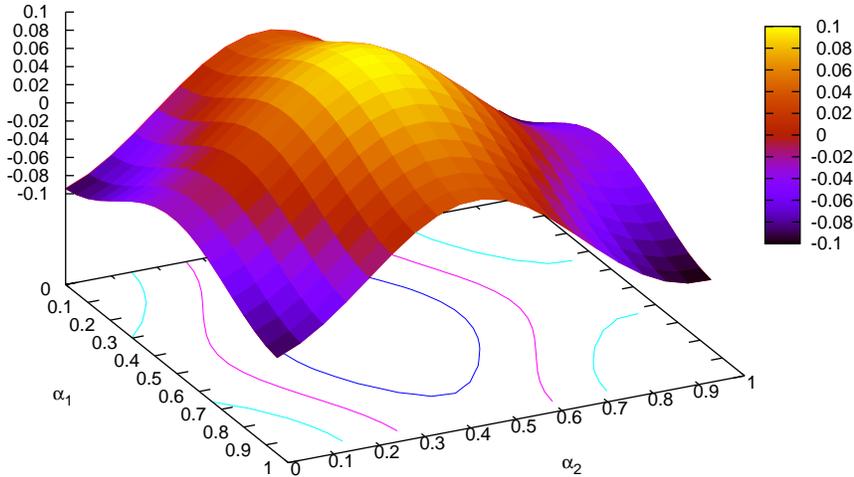}}
\caption{Jacobi Abel $\Psi$ function in the case of a set composed of three intervals.}
\label{fig-teta1}
\end{figure}

This result is the basis of our numerical approach to the isospectral torus.
For instance, Figure \ref{fig-teta1} displays the function $\Psi(\alpha_1,\alpha_2)$ in a case composed of three intervals, {\em i.e.} two fundamental frequencies, a case that permits a graphical rendering in three dimensions. How to parameterize conveniently this function is the topic of our further development.

Let us from now on focus on the sequence of the outdiagonal matrix Jacobi matrix elements $\{b_j\}_{j=1}^\infty$, that is more significant mathematically than the sister $\{a_j\}_{j=0}^\infty$. All results will nonetheless apply also to this latter.
Let us expand the analytic function $\Psi$ in Fourier series. Indicating by $\alpha$ the vector of arguments $\alpha_1,\ldots,\alpha_{N-1}$, a $N-1$ dimensional real vector, and by $k$ a $N-1$ dimensional integer vector, we can write
\begin{equation}
\label{nor3}
   \Phi (\mathbf{\alpha}) = \sum_{\mathbf{k} \in {\cal K}} C_\mathbf{k} \exp{ (2 \pi i \mathbf{k} \cdot \mathbf{\alpha})}.
\end{equation}
The set of indices ${\cal K}$ is to be chosen so to reflect the reality of the function $\Psi$.
Using this representation in eq. (\ref{eq:phi1}), one trivially finds that
\begin{equation}
\label{nor05}
   b_j  = \Psi (j \omega + \psi) =
   \sum_{\mathbf{k} \in {\cal K}} C_\mathbf{k}
   \exp{ (2 \pi i \mathbf{k} \cdot  \psi)}
   \exp{ (2 n \pi i \mathbf{k} \cdot \omega)}
   .
\end{equation}
However simple, the above equation reflects a significant, well known fact:
\begin{lemma}
\label{lem:alper1}
The outdiagonal sequence $\{b_j\}_{j=1}^\infty$ of a Jacobi matrix in the isospectral torus is almost periodic, with frequency module generated by the fundamental frequencies $\omega^n_i$.
\end{lemma}
In equation (\ref{nor05}), one observes the frequencies $\omega_\mathbf{k}:= 2 \pi  \mathbf{k} \cdot \omega$ and the phases $\psi_\mathbf{k} := 2 \pi \mathbf{k} \cdot  \psi$. Recall that the function $\Psi$ does {\em  not} depend on the specific matrix in the isospectral torus. Therefore, we discover that
\begin{lemma}
The almost--periodic spectrum of the sequence $\{b_j(\xi,\sigma)\}_{j \in \mathbf{N}}$  for any choice of $\xi$, $\sigma$ is characterized by the same set of frequencies;
the coefficient of the frequency $\omega_\mathbf{k}$ is $C_\mathbf{k} \exp{ (i \psi_\mathbf{k} )}$ and therefore it has a constant modulus, irrespective of the choice of $\xi$, $\sigma$.
\label{lem-freq}
\end{lemma}
The coefficient $C(\omega)$ of an almost periodic sequence $\{b_j\}_{j \in \mathbf{N}}$ is defined formally as the limit
\begin{equation}
  C(\omega) = \lim_{J \rightarrow \infty} \frac{1}{J} \sum_{j=1}^{J}
    e^{-i \omega j}  b_j.
  \label{eq:defin}
  \end{equation}

We can therefore restate our approach to the isospectral torus associated with the set $E^n$ as follows:
\begin{corollary}
The Jacobi-Abel function $\Psi$, represented as in eq. (\ref{nor3}) can be determined via a Fourier analysis of the sequence $\{b_j(\xi,\sigma)\}_{j=1}^\infty$, for any single set of values $\xi,\sigma$.
\end{corollary}

To be rendered effective, this corollary needs a refined Fourier analysis of the sequence \bseq. This can be achieved in three steps.
\begin{itemize}
\item[$\bullet$] We compute the frequencies $\omega_\mathbf{k}$. For this, we solve the set of equations (\ref{meas5}) and we compute the fundamental frequencies $\omega^n_i$, following \cite{dolo,arx}.
\item[$\bullet$] Using this information, we consider the numerical problem of computing the coefficients $C(\omega_{\bf k})$ in a more effective way than by the definition (\ref{eq:defin}) that is theoretically significant but numerically ineffective. A technique to this aim is described in Appendix, when the set of indices ${\cal K}$ is replaced by a finite set.
   \item[$\bullet$] Finally, we organize the infinite set of indices as an increasing sequence of finite, nested sets.
   \end{itemize}

\begin{figure}
\centerline{\includegraphics[width=.6\textwidth, angle = -90]{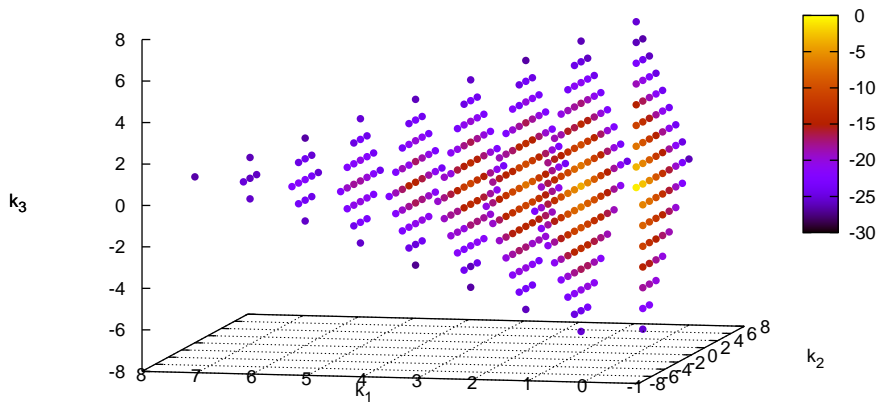}}
\caption{Amplitudes $|C_\mathbf{k}|$ for $E^2$, example \ref{examp-cantor}.}
\label{fig-ampli}
\end{figure}

The first two steps above being described in Appendix and in the references, let us concentrate on the third. A possibility in this respect is to use a metric in ${\mathbf Z}^{N-1}$, that includes ${\cal K}$, so to define
\[
{\cal K}_L = {\cal K} \cap B_L
\]
where $B_L$ is the ball of radius $L$ centered at the origin. Convergence of the Fourier sequence (\ref{nor3}), when ${\cal K}$ is replaced by
${\cal K}_L$ and $L$ tends to infinity is a classical problem of harmonic analysis, see {\em e.g.} \cite{feffer}. In Figure \ref{fig-ampli} we use the metric induced by the norm $\| \mathbf{k} \| = \sum_{i=1}^{N-1} |k_i|$, we set $L=8$ and we consider $E^n$, with $n=2$. In this case the support of the measure is made of four intervals: there are three fundamental frequencies $\omega^2_i$ and $\mathbf{k}$ is a three--dimensional vector. At the location of this vector in space, we display a color-coded point, in logarithmic scale, according to the amplitude $|C_\mathbf{k}|$.
We observe an exponential decrease of the amplitudes as $\|\mathbf{k}\|$ increases, that supports the relevance of the metric adopted. Yet, at an accurate inspection, this decrease is not uniform, but this observation leads to further discovery.

\begin{figure}
\centerline{\includegraphics[width=.6\textwidth, angle = -90]{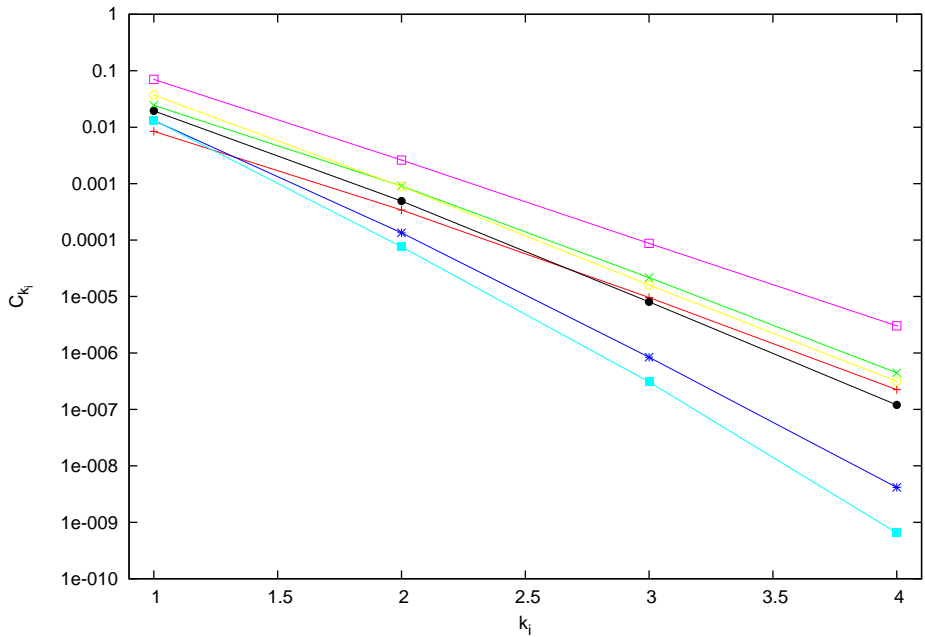}}
\caption{Amplitudes $|C_\mathbf{k}|$ for $E^3$, example \ref{examp-cantor}, versus $k_i$, when $\mathbf{k} = (0,\ldots,k_i,0,\ldots)$, $i =1,\ldots,7$.}
\label{fig-ampli2}
\end{figure}

In fact, let us plot the values of $|C_\mathbf{k}|$ along lines from the origin directed as the coordinate axes. That is, fix all $k_i$ values to zero except one, and raise this latter from $1$ to $L$. We obtain fig. (\ref{fig-ampli2}), now for $E^3$. The behavior observed for $n=3$ is typical, and consists of an exponential decay of amplitudes with $k_i$.

\begin{figure}
\centerline{\includegraphics[width=.6\textwidth, angle = -90]{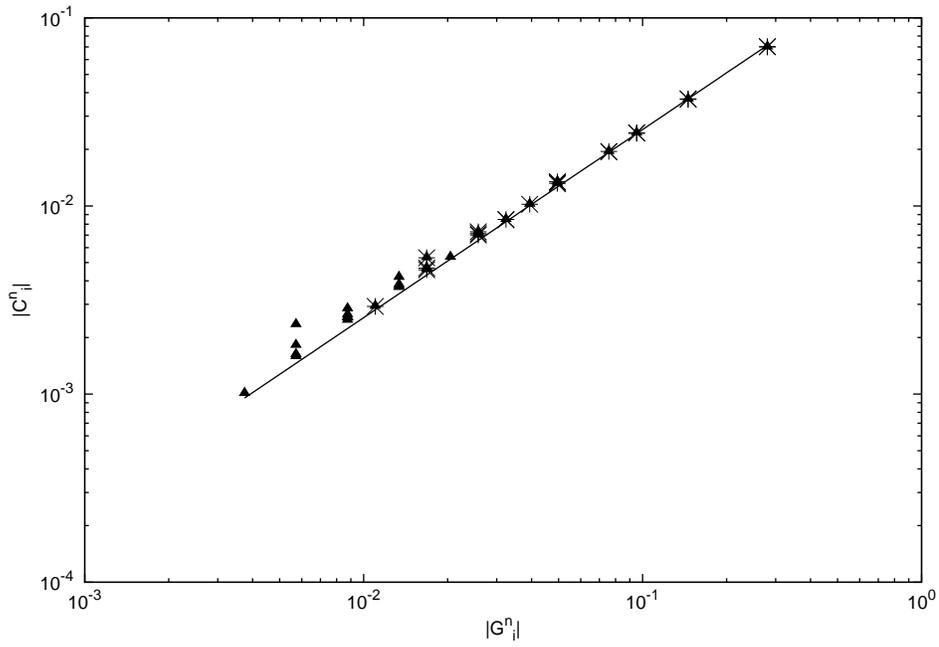}}
\caption{Principal harmonic amplitudes $|C^n_i|$ versus $|G^n_i|$, the size of the $i$-th gap in $E^n$. Values for $n=2,\ldots,5$ are indicated by different symbols, example \ref{examp-cantor}. The fitting line is $g(x) = A x$, with a suitable constant $A$. }
\label{fig-gappa1}
\end{figure}

The decays starts off from the value of the amplitude of the principal harmonic $\omega^n_i$, that for simplicity we denote $|C^n_i|$. This harmonic is associated with the gap $G^n_i$ in the set $E^n$. It has been observed long ago in the theory of electronic states in solids, that an harmonic component in the potential of a tight binding Hamiltonian opens a gap in the spectrum of this latter. By perturbation theory arguments \cite{duca} it was established that the size of the gap is proportional to the amplitude of the harmonic component. We observe here the reverse behavior: proportionality is approximatively well verified, but it is the size of the gap that induces the amplitude of the principal harmonic: see fig. \ref{fig-gappa1}.

Remark that data points in figure \ref{fig-gappa1} refer to four different orders $n$. Also notice that some points at different values of $n$ appear to be almost coincident. These observations are of relevance when considering the limit of infinite $n$, a thing that we shall do later on.

\begin{figure}
\centerline{\includegraphics[width=.6\textwidth, angle = -90]{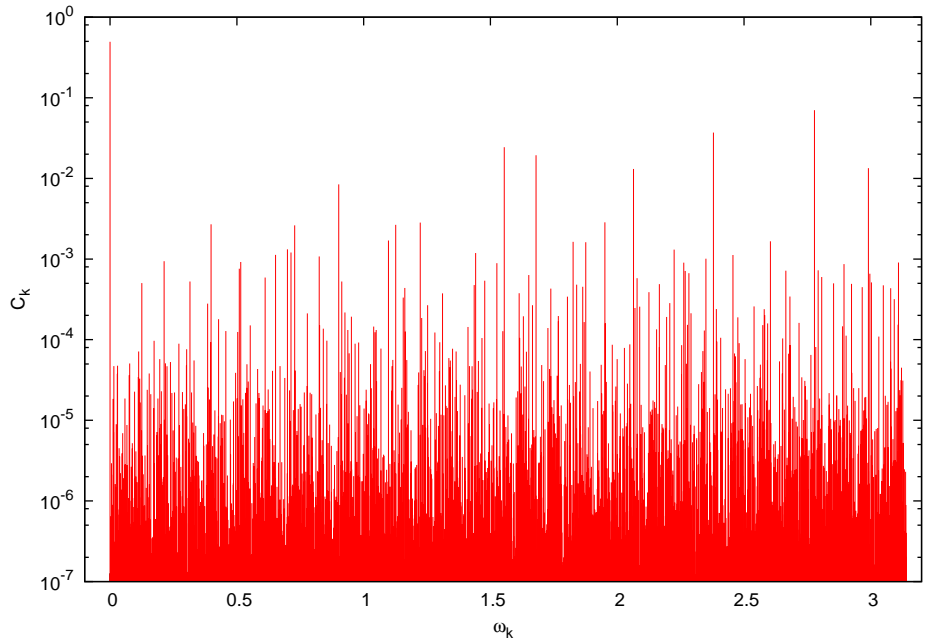}}
\caption{Leading harmonic amplitudes $C_{\mathbf k}$ versus $\omega_{\mathbf k}$, for  $n=3$, $\|{\mathbf k}\| \leq 6$.}
\label{fig-ampliz}
\end{figure}

To conclude this section, we show in fig. \ref{fig-ampliz} the harmonic amplitudes
$C_{\mathbf k}$ versus $\omega_{\mathbf k}$ obtained by our technique for $n=3$.
Remark two fundamental differences with respect to fig. 2 in \cite{physd1} (or the analogous fig. 3.1 in \cite{simonlast}): a difference in technique---fast Fourier was there employed, limiting resolution, and a quite evident difference in symmetry.
We are now in the presence of a typical spectrum (although still at a low value of $n$) of an IFS system.  Also to test the reliability of the data so displayed, in the next section we want to tackle a difficult problem, that of convergence to the isospectral torus.

\section{Convergence to the isospectral torus}
\label{sec-harm}

In this section we will study the convergence of Jacobi matrices of different measures supported on $E^n$ to their limit point in the isospectral torus. We will consider the following situation:
\begin{hypothesis}
Suppose that $\mu_n$ is an absolutely continuous measure with respect to the Lebesgue measure supported on $E^n$. Also suppose that $\mu$ satisfies a Szeg\"o condition on $E^n$.
\label{hyp-sze}
\end{hypothesis}
We denote by  $\{a_j(\mu_n)\}_{j \in \mathbf N}$,  $\{b_j(\mu_n)\}_{j \in \mathbf N}$ the sequences of Jacobi matrix elements of $\mu_n$. In this respect, the following theorem is notable

\begin{theorem}[\cite{franzinfty}]
Under hypothesis \ref{hyp-sze} the recurrence coefficients $a_j(\mu_n)$ and $b_j(\mu_n)$ are asymptotically convergent to the almost periodic form given by eq. (\ref{eq:phi1}), where the functions $\Phi$, $\Psi$ and the vector $\omega$ are derived from $E^n$, while the shift vectors $\phi$ and $\psi$ depend on $\mu$.
\label{teofranz}
\end{theorem}

A comment is in order. The above theorem has been proven in \cite{franzinfty} in a much more general setting than the one presented here. On the other hand, a strengthening of this result is possible for special classes of measures. Aptekarev \cite{sasha} considered absolutely continuous measures on $E^n$ with densities  on each interval $E^n_i$ of the kind $\rho(s) = w(s) (s-a_i)^{\alpha_i}(b_i-s)^{\beta_i}$, being $\alpha_i,\beta_i>-1$ and $w(s)$ locally analytical and non-negative. He proved that distance from the asymptotic limit provided by eq. (\ref{eq:phi1}) is asymptotically $O(1/j)$. Suetin \cite{suetrace}, by requiring that $\rho(s) = w(s)/\sqrt{|Y(s)|}$ with a strictly positive, holomorphic $w$, proves that convergence is $O(\delta^n)$, with a suitable constant $\delta<1$.
We will now study numerically both cases.

To define measures satisfying Aptekarev or Suetin conditions, we make recourse once more to the IFS construction. Let us associate to each map $\phi_j$ in eq. (\ref{mappi}) a positive weight $\pi_j$, $\sum_{j=1}^M \pi_j = 1$ and let us define
the operator $T$ acting on continuous functions on $\mathbf R$ given by $T f  = \sum_j \pi_j f \circ \phi_j$. Let $T^*$ be the dual operator of $T$, acting in the space of regular Borel probability measures. Then, the following theorem is standard
\begin{theorem}[\cite{dem}]
Let $\mu_0$ any Borel probability measure on $E^0=[\min \gamma_{j},\max \gamma_j]$. Let $\mu_n = (T^*)^n \mu_0$. Then $\mu_n$ tends weakly and in Kantorovich distance to a measure $\mu$ supported on the attractor ${\cal A}$, that is called an {\em IFS balanced} measure and satisfies the fixed point condition $T^* \mu = \mu$.
\label{teo-meas}
\end{theorem}

The devil's staircase/Cantor measure considered in the Introduction is precisely of this kind. Observe that Theorem \ref{teo-meas} applies to {\em any} initial measure $\mu_0$.
We therefore make two choices. In case a) we use the normalized Lebesgue measure on $E^0$. In case b) $\mu_0$ is the (suitably rescaled) orthogonality measure of  Chebyshev polynomials of second kind.
It then easily follows that
\begin{lemma}
Let $\mu_n$ be defined as above, with the IFS maps (\ref{mappi}) and initial measure $\mu_0$ as in cases a), b). Then, for any $n$, $\mu_n$ is an absolutely continuous measure that verifies hypothesis \ref{hyp-sze}.
In case a) and b) Aptekarev condition is satisfied, while Suetin condition holds only for case b).
\end{lemma}

Let us now apply a numerical technique for the stable computation of the Jacobi matrices $J_{\mu_n}$, that applies to both cases a) and b) just discussed.
This technique has been introduced in ref. \cite{mobius} and is radically different from the one described in \cite{arx} and in Sect. \ref{sec-equi} of this work.
On the force of the above theorems, these matrices converge to a point in the isospectral torus. This raises a challenging numerical problem.

\begin{problem}
Given a measure $\mu$ satisfying hypothesis \ref{hyp-sze}, and its Jacobi matrix $J_\mu$, find the matrix in the isospectral torus of $E^n$ to which $J_\mu$ converges.
\label{prob-1}
\end{problem}

We propose the following solution of this problem.
In view of the harmonic analysis presented in Sect. \ref{harman} and of Lemma \ref{lem-freq}, we can reduce the previous problem to the determination of the shifts $\psi_{\mathbf k}$, for $\mathbf{k} \in {\cal K}_L$. It is then possible from these latter to reconstruct the precise isospectral matrix, since amplitudes and frequencies are universal in the isospectral torus and can be determined in advance. We denote by $\theta_n$ and $J_{\theta_n}$ measure and matrix in the isospectral torus that correspond to $\mu_n$.

Again, let us consider for simplicity the single sequence $\{b_j(\mu_n)\}_{j \in \mathbf{N}}$. The associated set of phases $\psi_{\mathbf k}$ can be computed combining Fourier analysis with a suitable extrapolation procedure outlined in the second part of the Appendix. Let us comment here on the results that can be  obtained in this way.

\begin{figure}
\centerline{\includegraphics[width=.6\textwidth, angle = -90]{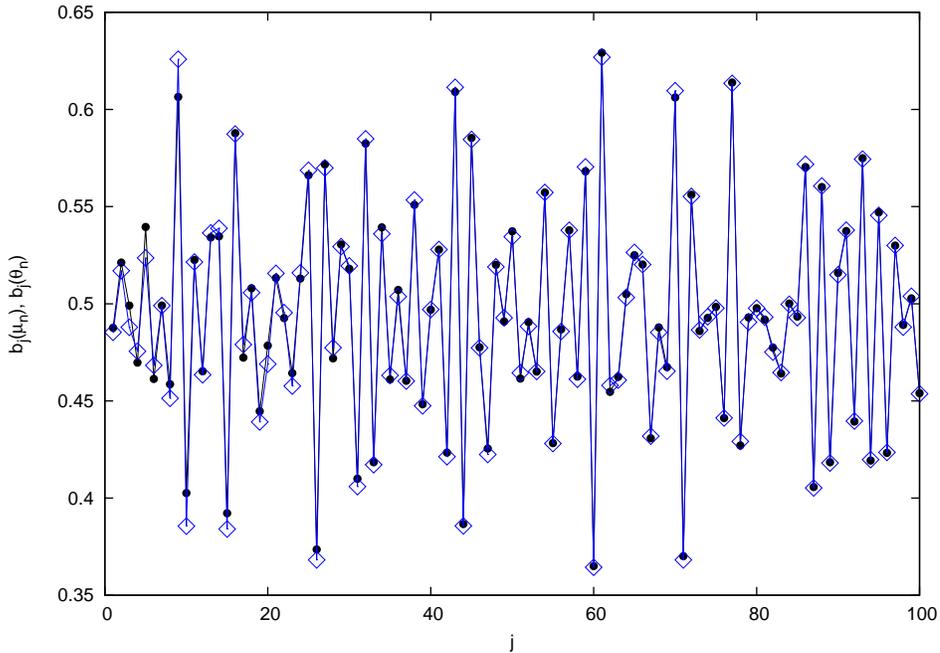}}
\caption{Outdiagonal entries of the Jacobi matrix of $\mu_n$ (circles) and of its limit $\theta_n$ in the isospectral torus (diamonds). Here, $n=4$ and we have used
example \ref{examp-cantor} with $\pi_1=3/5, \pi_2=2/5$ and initial measure $\mu_0$ in case a).}
\label{fig-convepoint}
\end{figure}

We start by examining case a). Fig. \ref{fig-convepoint} exhibits convergence of $b_j(\mu_n)$ to the corresponding matrix element $b_j(\theta_n)$ of the measure $\theta_n$ in the isospectral torus of $E^n$, obtained as described above. Convergence is slow enough to permit us to appreciate visually how the two sequences are getting close.

\begin{figure}
\centerline{\includegraphics[width=.6\textwidth, angle = -90]{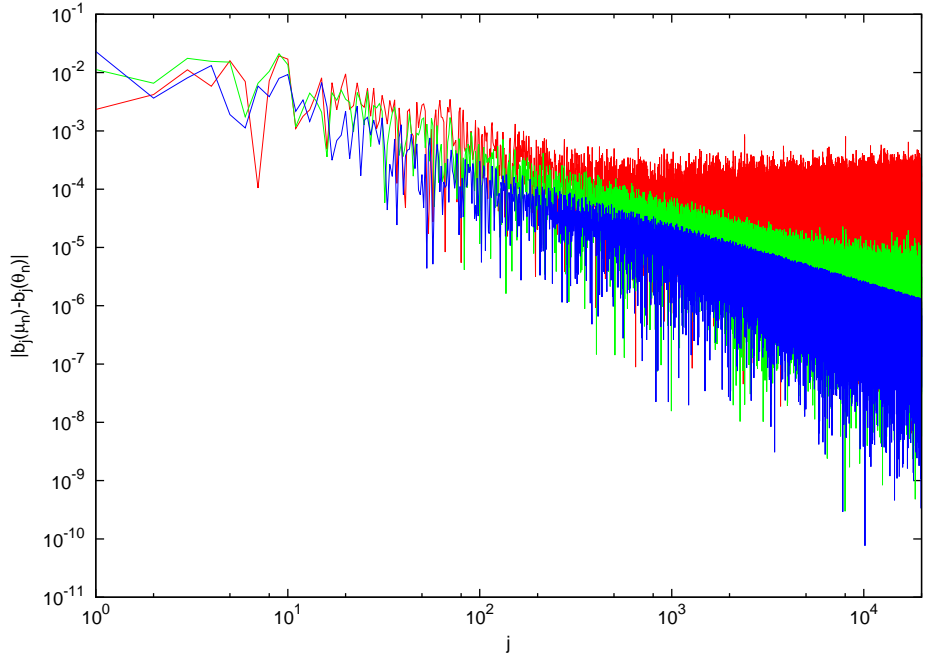}}
\caption{Difference between the outdiagonal entries of the Jacobi matrix of $\mu_n$ and of its limit $\theta_n$ in the isospectral torus. Here, $n=2$ (blue line), $n=3$ (green) and $n=4$ (red). We have used
example \ref{examp-cantor} and $\mu_0$ in case a).}
\label{fig-manydec1}
\end{figure}

To quantify convergence, we compute the difference $|b_j(\mu_n)-b_j(\theta_n)|$ and we plot it in double logarithmic scale in fig. \ref{fig-manydec1}. We observe that the decay of this difference is bounded from above by $A_n/j$, in accordance with \cite{sasha}. At large values of $j$ this decay stops and difference between $b_j(\mu_n)$ and $b_j(\theta_n)$ oscillates at a low value due to the combination of two effects: on the one hand, computations are performed with a finite set of frequencies ${\cal K}_L$. When $n$ increases the cardinality of ${\cal K}_L$ increases as $L^n$, so that we are forced to reduce $L$ in order to keep the number of frequencies within acceptable limits (see Appendix). On the other hand, our Fourier analysis is performed on finite Jacobi matrices whose length is fixed, albeit large (typically, 800,000). Therefore, amplitudes and phases of $\omega^n_{\mathbf k}$ are determined less precisely when $n$ grows.

\begin{figure}
\centerline{\includegraphics[width=.6\textwidth, angle = -90]{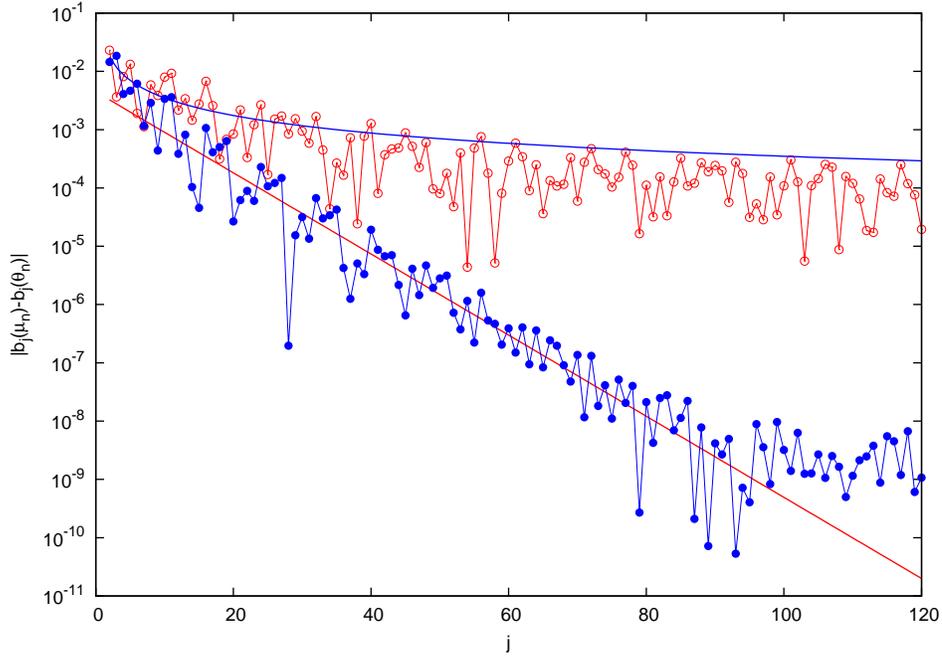}}
\caption{Difference between the outdiagonal entries of the Jacobi matrix of $\mu_n$ and of its limit $\theta_n$ in the isospectral torus. Here, $n=2$ and $\mu_0$ is case a) (open red circles) and b) (filled blue circles). The curves plotted as reference (in reverse color coding to improve visibility) are an exponential (red) and a power-law decay $A/j$ (blue).}
\label{fig-deche2-1}
\end{figure}

\begin{figure}
\centerline{\includegraphics[width=.6\textwidth, angle = -90]{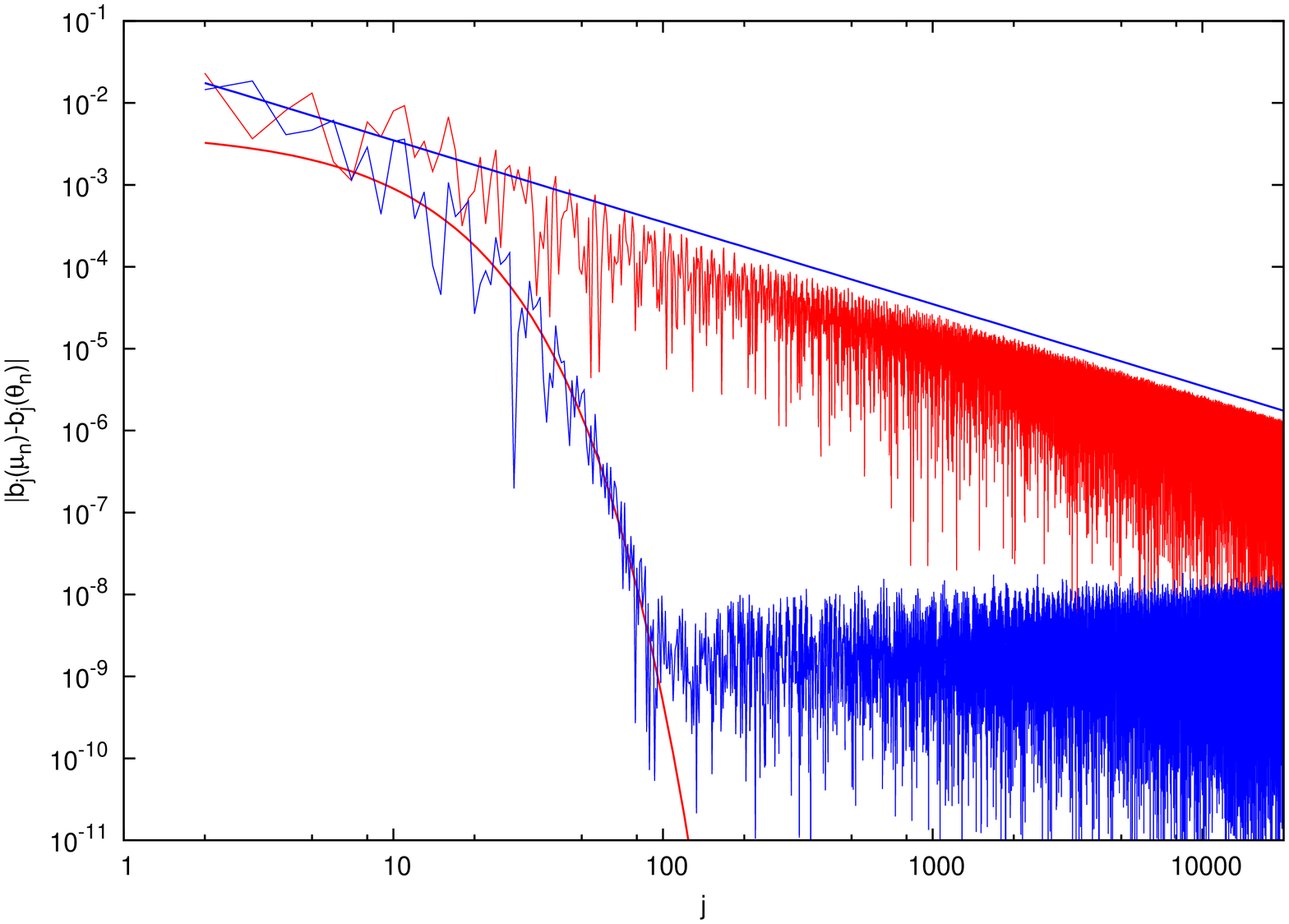}}
\caption{Same data as in fig. \ref{fig-deche2-1}, now in double logarithmic scale.}
\label{fig-deche2-2}
\end{figure}

Next, let us consider case b). Figure \ref{fig-deche2-1}  is the analogue of fig. \ref{fig-manydec1}, but it now compares cases a) and b) at the fixed value of $n=2$.
The exponential decay of the latter is immediately appreciated.
Plotting the same data on a larger interval, and resorting to logarithmic scale, shows the $O(1/j)$ bound on the decay: see fig. \ref{fig-deche2-2}. The exponential convergence of case b) appears even more dramatically here. The same comments on precision of the computation described above also apply in this case. 

We can now conclude this section noting that the data presented confirm the validity of our approach in determining numerically the isospectral torus for finite collections of intervals. We can therefore move to the Cantor set case.

\section{The isospectral torus of infinite gap sets--Cantor sets}
\label{sec-isoinfty}
A recent paper by Kr\"uger and Simon \cite{simonlast} discusses in the usual illuminating way various issues related to the definition(s) of an isospectral ``torus'' associated with measures supported on Cantor sets. Our approach to this problem is the following: in accordance with definition \ref{defi-isofinite} we would like to define this object as
\begin{definition}
The isospectral torus of a hyperbolic IFS Cantor set with finitely many maps is the set of Jacobi matrices associated with measures $\theta_\infty$, defined as weak limits of $\theta_n$, when this latter is defined in eqs. (\ref{lupa1},\ref{lupa2},\ref{meas6}) and $E^n$ is defined in eq. (\ref{eq-int1}).
\label{defi-isoinfinite}
\end{definition}

It is clear that the above definition can be extended beyond the realm of I.F.S. with a finite number of maps. Yet, we have been intentionally conservative, postponing generalizations. In fact, we believe that these hypotheses should be enough to assure weak convergence of the following procedure.
Write the gaps in $E^n$ as:
 \begin{equation}
\label{eq-int2}
 G^n = G^{n-1} \bigcup H^n,
\end{equation}
where $H^n$ is the set of ``new'' gaps created at level $n$. $H^n$ is iteratively constructed as
\begin{equation}
\label{eq-gap2}
  H^{n} = \bigcup_{j=1}^{M}  \phi_j(H^{n-1}),
\end{equation}
starting from $H^1 = (\beta^1_1,\alpha^1_2) \bigcup \ldots \bigcup (\beta^1_{M-1},\alpha^1_M)$.
The countable set of gaps $G^\infty$ in the I.F.S. Cantor set is then written as
\begin{equation}
\label{eq-gapinf}
  G^{\infty} = \bigcup_{n=1}^{\infty} H^{n}.
\end{equation}
In so doing, we can order the countable set of gaps of the Cantor set as they appear with $n$ in eq. (\ref{eq-gapinf}), and for each $H^n$ in ascending order on the real line. This also imply an ordering of the infinite set of points $\{\xi_i\}_{i \in \mathbf N}$ (one in each gap) and of the Riemann surface indices $\{\sigma_i\}_{i \in \mathbf N}$.

Let us now consider such a double sequence $\{\xi_i,\sigma_i\}_{i \in \mathbf N}$, that defines a sequence of measures $\{\theta_n\}_{i \in \mathbf N}$. If this sequence converges weakly, then the Jacobi matrix elements $a_j(\theta_n)$ and $b_j(\theta_n)$ also converge. As mentioned, we conjecture that IFS hyperbolicity should warrant this convergence.

\begin{figure}
\centerline{\includegraphics[width=.6\textwidth, angle = -90]{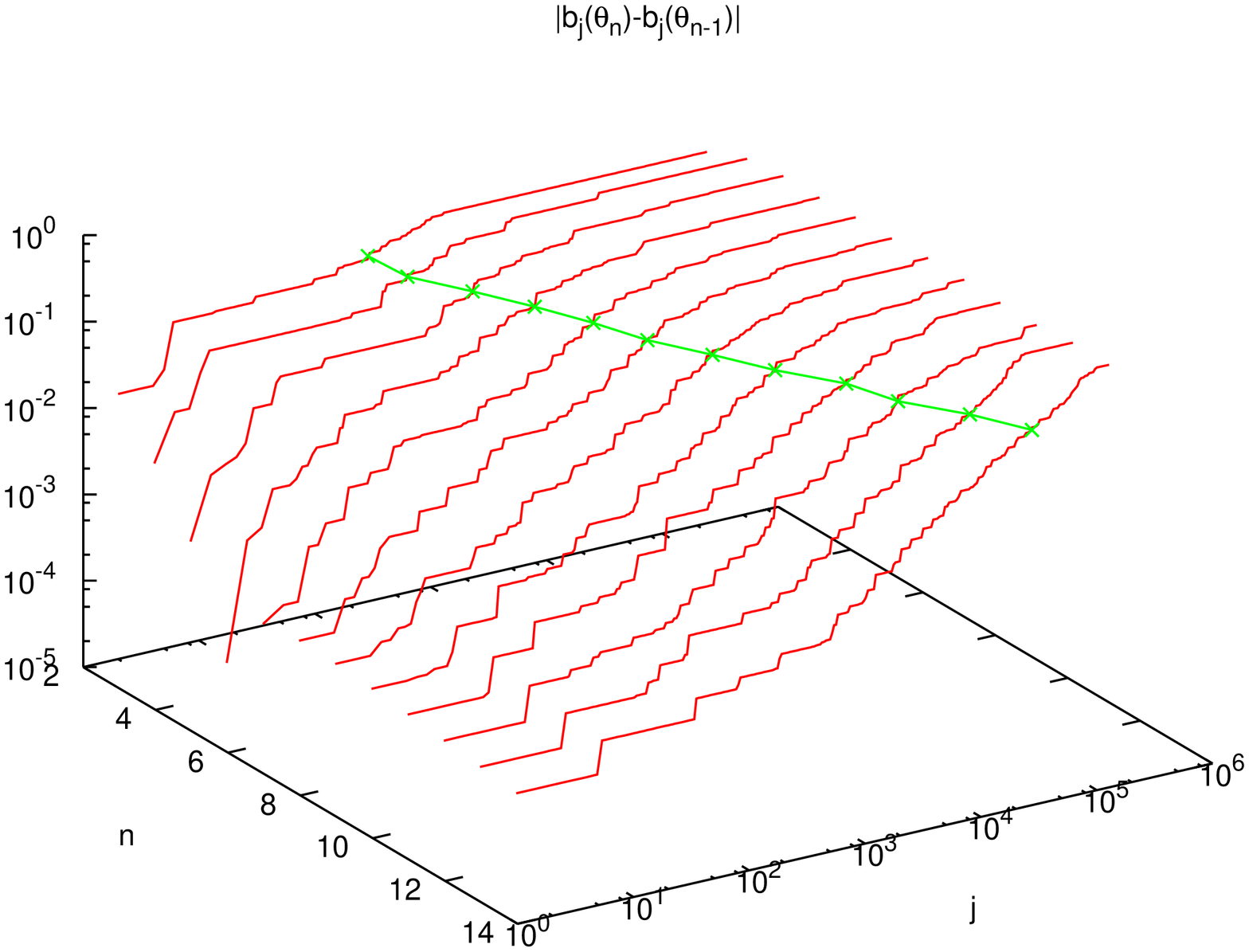}}
\caption{Computing the isospectral torus of an infinite--gap set: see text for details.}
\label{fig-torconve}
\end{figure}

Our technique permits to test this conjecture numerically. As a matter of fact,  convergence seems to hold at a geometrical rate---as it might be expected, because of  the hyperbolicity hypothesis. In fig. \ref{fig-torconve} we plot the absolute value of the difference $b_j(\theta_n)-b_j(\theta_{n-1})$, versus both $n$ and $j$. On the surface defined by these values, we also plot the line $N(\varepsilon,n)$ indicating that for $j \leq N(\varepsilon,n)$ the difference between $b_j(\theta_n)$ and its infinite limit is less than $\varepsilon$. This quantity is significant, because it indicates that, {\em at precision $\varepsilon$, the truncated isospectral Jacobi matrix of rank $N(\varepsilon,n)$ is insensitive to the torus variables $\{\xi_i,\sigma_i\}_{i > Q}$}, with $Q=M^n$.
As a consequence of geometric convergence, for any $\varepsilon$, $N(\varepsilon,n)$ grows asymptotically in an exponential way with $n$.

The same geometric convergence is observed in the harmonic analysis of these Jacobi matrices. This holds for both frequencies and amplitudes. For example, in Table \ref{tab-1} we report the fundamental frequency corresponding to the second gap in our ordering, and its amplitude. This geometrical convergence might be the key by which to resolve the conjecture mentioned in the introduction, at least numerically. We now turn to this last point in this paper.

\begin{table}
\centering
\begin{tabular}{|c|l|l|l|l|}
  \hline
   & $n=2$ & $n=3$ & $n=4$ & $n=5$  \\
  \hline
  $\omega^n_i$  &  1.55434055   & 1.55445514   &  1.55438155  &    1.55429495 \\
   $C^n_i$ &  2.42975661E-002   &   2.44107202E-002   & 2.44958535E-002     &  2.45456487E-002    \\
  \hline
\end{tabular}
\caption{Fundamental frequency $\omega^n_i$ and related amplitude $C^n_i$ for the second gap in $E^n$, for the IFS of example \ref{examp-cantor}. \label{tab-1}}
\end{table}

\section{If they are almost--periodic?}
\label{sec-if}\footnote{The title of this section was clearly derived from \cite{belli2}.}
As mentioned, the numerical techniques presented in this paper aim at providing hints to the resolution of the conjecture on the almost periodicity of hyperbolic IFS Jacobi matrices. We want to report here preliminary results in this direction. These results are {\em not} resolutive in a sense or another, but clarify the setting of the problem.

For sake of definiteness, let $\mu_\infty$ be the balanced measure of such an IFS, $J_\mu$ its Jacobi matrix, and $\{b_j(\mu)\}_{j \in \mathbf N}$ the sequence of its out-diagonal entries. Let also $\mu_n$ be the sequence of its IFS approximations, constructed starting from the Chebyshev measure of second kind, and
$\theta_n$ the corresponding measures in the {\em finite dimensional} isospectral torus.

\begin{figure}
\centerline{\includegraphics[width=.6\textwidth, angle = -90]{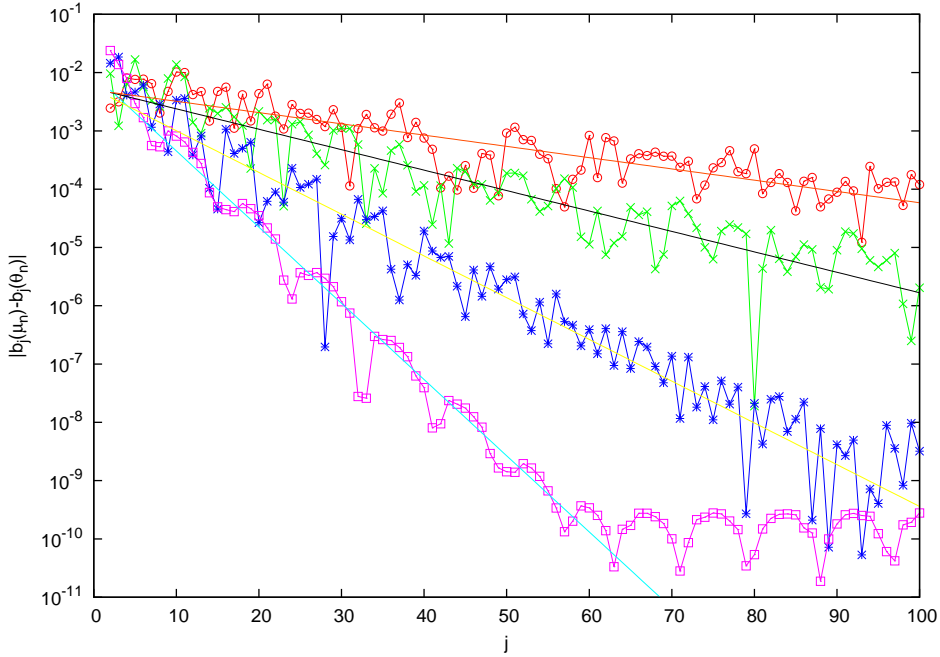}}
\caption{Differences $|b_j(\theta_n)- b_j(\mu_n)|$ versus $j$, for $n=1,\ldots,4$, and fitting exponential curves. }
\label{fig-slopeg1b}
\end{figure}

We have seen that as $j$ grows, $b_j(\theta_n)$ and $b_j(\mu_n)$ approach each other exponentially fast. It is now important to gauge the {\em slope} of this exponential. In figure \ref{fig-slopeg1b} we plot the cases $n=1$ to $n=4$ for the example considered in this paper. The fitting exponential lines have form $c_n e^{-d_n j}$. They have been computed as a log-linear least squares fit, and appear to be approximately intersecting at $n \sim 1$. The important feature is then the decay coefficient $d_n$. In table \ref{tab-2} we report these values. It appears that, as $n$ grows, $d_n$ tends to zero. As a matter of fact, the sequence of $d_n$ values can be itself extremely well approximated by an exponential decay: $d_n \sim e^{-\delta n}$ with $\delta \simeq 0.644$.

Let us now fix a positive threshold $\varepsilon$. Extrapolating the data so far obtained, we can predict that for $j > J(\varepsilon,n) = \log(1/\varepsilon) e^{\delta n}$ the distance between $b_j(\theta_n)$ and $b_j(\mu_n)$ is less than $\varepsilon$.

\begin{table}
\centering
\begin{tabular}{|c|l|l|l|l|}
  \hline
   & $n=1$ & $n=2$ & $n=3$ & $n=4$  \\
  \hline
  $d_n$  &  0.300994     & 0.164903     &  0.080720    &    0.044511  \\
  \hline
\end{tabular}
\caption{Coefficients of the fitting exponentials in fig. \ref{fig-slopeg1b}}. \label{tab-2}
\end{table}

Let us now consider the distance between $b_j(\mu_n)$ and $b_j(\mu_\infty)$. At difference with the previous, this distance can be bound by a fixed value $\varepsilon$ when $j$ is {\em smaller} than a value $N(\varepsilon,n)$.
We therefore define $N(\varepsilon,n) := \max \{l \mbox{  s.t. } |b_j(\mu_n)-b_j(\mu_\infty)| \leq \varepsilon,  0 < j \leq l \}$.
It was originally found in \cite{mobius} that $N(\varepsilon,n)$ increases exponentially with $n$. Notice that the same dependence has been observed in Sect. \ref{sec-isoinfty} and fig. \ref{fig-torconve}, but we are now in a completely different setting, because one of the terms, $\mu_\infty$, is the I.F.S. singular continuous measure that is the object of the debated conjecture.
Suppose that it happens that for $\varepsilon$ sufficiently small, $N(\varepsilon,n) \sim e^{\kappa n}$. Then, if the inequality $\kappa > \delta$ holds, by choosing larger and larger values of $n$, the Jacobi matrix elements $b_j(\mu_\infty)$ can be uniformly approximated by $b_j(\theta_n)$ in the asymptotically long, asymptotically remote segments $j \in (e^{\delta n},e^{\kappa n})$.

\begin{figure}
\centerline{\includegraphics[width=.6\textwidth, angle = -90]{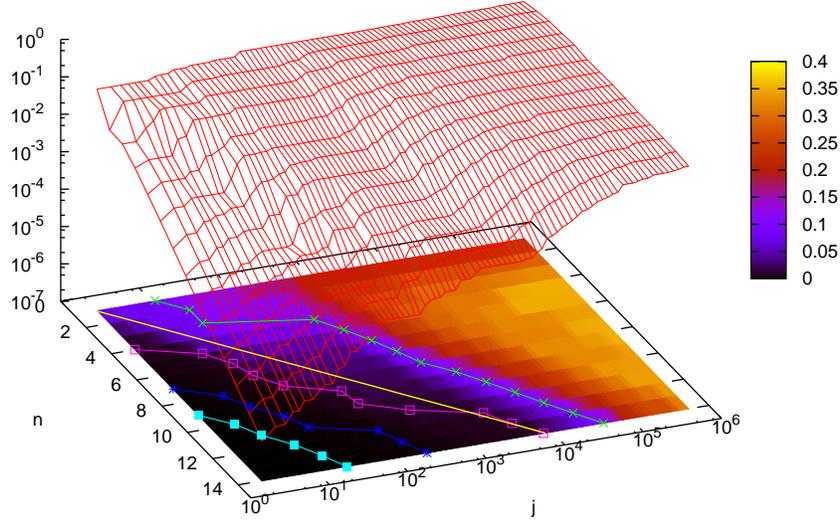}}
\caption{Differences $|b_j(\mu_n)- b_j(\mu_\infty)|$ versus $n$ and $j$. Also plotted at the basis of the picture are the values $j(n)=N(\varepsilon,n)$, for $\varepsilon=10^{-1},10^{-2},10^{-3},10^{-4}$, and the (yellow) line $k(n)=e^{\delta n}$, with $\delta=0.644$.}
\label{fig-final2b}
\end{figure}

Our original technique \cite{cap} of computing the Jacobi matrix associated with $\mu_\infty$ can be now put to work, in combination with that in \cite{mobius}, to compute the function $|b_j(\mu_n)-b_j(\mu_\infty)|$, that is plotted logarithmic scale in fig. \ref{fig-final2b}.  At the basis, the experimental curves $j(n)=N(\varepsilon,n)$ for various values of $\varepsilon$ are plotted, as well as $k(n)=e^{\delta n}$, with $\delta$ given by the extrapolation of the data in table \ref{tab-2}. Even if the figure spans over various orders $n$ and up to Jacobi matrices of size in the hundreds of thousands, a clear-cut evidence in favor of the inequality $\kappa>\delta$, or its violation, cannot be obtained: the displayed data seem to be very close to the hedge between the two behaviors.

\section{Appendix: Discrete Fourier Analysis}

The harmonic problem required in this paper differs from the usual in that the frequencies present in the signal are known, being determined by eq. (\ref{nor05}).
Therefore, we act as follows. First, we multiply the ``signal'' $b_j$ by a window function $w_j$. For this, we choose the Dolph--Chebyshev window, well described in the classical paper \cite{harris}. Define the windowed discrete Fourier transform
\begin{equation}
  F(\omega) = \sum_j b_j w_j e^{-i \omega j}.
 \label{cheby1}
 \end{equation}
Due to the property of the filter, the above is a finite summation.
Next, we take into account the fact that the frequencies $\omega_{\mathbf k}$ are known.
Since we restrict ${\mathbf k}$ to ${\cal K}_L$, this set is finite, and we can list frequencies in ascending order. Evaluating the windowed Fourier transform at any of these frequencies we find that
\begin{equation}
  F(\omega_{\mathbf k}) = \sum_{{\mathbf k'} \in {\cal K}_L} W({\mathbf k},{\mathbf k'}) D_{\mathbf k'},
 \label{cheby2}
 \end{equation}
where $D_{\mathbf k'}$ is a complex coefficient that encodes both $C_{\mathbf k'}$ and $\psi_{\mathbf k'}$. Actually, we work with two separate equations for real and imaginary part: this has several numerical advantages.

$W({\mathbf k},{\mathbf k'})$ is a kernel that can be exactly computed, so that eq. (\ref{cheby2}) is a linear system for the unknowns $D_{\mathbf k'}$. Taken at face value, this equation is discouraging, for the cardinality of the set ${\cal K}_L$  makes it untractable numerically. Yet, the Dolph-Chebyshev filter is such that $W$ is a function of the difference $\omega_{\mathbf k}-\omega_{\mathbf k'}$, and this function can be made arbitrarily small outside of an interval centered at zero \cite{harris}. Therefore, in our numerical computations, by arranging frequencies in increasing order, $W({\mathbf k},{\mathbf k'})$ becomes a banded matrix, and the system (\ref{cheby2}) can be swiftly solved. It is just the case to remark that this procedure yields results of a precision unattainable by the mere usage of the fast Fourier transform, that was applied to this problem in \cite{physd1,simonlast}.

Finally, a remark on the problem of extracting the phases $\psi_{\mathbf k}$ when examining the sequence $b_j(\mu_n)$. In this case, the finite size window $w_j$ is shifted along the sequence, by a finite lag $l$. A sequence of lags $l_m$, $m=1,\ldots,M$ is performed, and the corresponding amplitudes and phases extracted. A fitting/extrapolation procedure is then performed that yields the asymptotic values $C_{\mathbf k}$ and $\psi_{\mathbf k}$. The former are then compared with those obtained by the isospectral torus to assess the validity of the procedure, while the latter yield the desired, unknown phases.


\end{document}